\documentclass{amsart}
\usepackage{amssymb,euscript,amsmath, mathrsfs}
\usepackage[dvips]{graphicx}
\usepackage[dvips]{color}
\usepackage{epsfig}

\newcounter{ENUM}
\newcommand{\itm}{\item}
\newenvironment{ilist}{\renewcommand{\theENUM}{\roman{ENUM}}\renewcommand{\itm}{\addtocounter{ENUM}{1}\item[(\theENUM)]}\begin{itemize}\setcounter{ENUM}{0}}{\end{itemize}}
\newenvironment{alist}[1][0]{\renewcommand{\theENUM}{\alph{ENUM}}\renewcommand{\itm}{\addtocounter{ENUM}{1}\item[\theENUM)]}\begin{itemize}\setcounter{ENUM}{#1}}{\end{itemize}}

\newcommand{\nHat}[1]{\widehat{#1}}                        

\def\1{{\mathbf 1}}

\def\EH{{\operatorname{EH}}}
\def\D{{\mathscr D}}

\def\j{{\mathfrak j}}

\def\N{{\mathbb N}}

\def\EH{\operatorname{EH}}

\newtheorem{thm}{Theorem}[section]
\newtheorem{prop}[thm]{Proposition}
\newtheorem{lem}[thm]{Lemma}
\newtheorem{cor}[thm]{Corollary}

\theoremstyle{definition}
\newtheorem{defn}[thm]{Definition}

\newtheorem{ex}[thm]{Example}

\theoremstyle{remark}

\newtheorem{rem}[thm]{Remark}

\numberwithin{equation}{section}

\address{Department of Mathematics, One Shields Avenue, University of California, Davis 95616.}
\email{fuliu@math.ucdavis.edu}

\begin{document}
\title{Moduli of Crude Limit Linear Series}
\author{Fu Liu}
\begin{abstract}
Eisenbud and Harris introduced the theory of limit linear series, and constructed a space parametrizing their limit linear series. Recently, Osserman introduced a new space which compactifies the Eisenbud-Harris construction. In the Eisenbud-Harris space, the set of refined limit linear series is always dense on a general reducible curve. Osserman asks when the same is true for his space. In this paper, we answer his question by characterizing the situations when the crude limit linear series contain an open subset of his space.
\end{abstract}

\maketitle

\section{Introduction}
Fix $r, d \in \N.$ In the 1980's Eisenbud and Harris introduced the concept of limit linear series of dimension $r$ and degree $d$ on a reducible curve $X$ of compact type and genus $g$ \cite{eisenbudharris1983}. 
They also constructed a space $G_d^{r, \EH}(X)$ parametrizing their limit linear series on $X.$ There is an open subset of $G_d^{r, \EH}(X)$ consisting of ``refined'' limit linear series, and its complement consists of ``crude'' limit linear series. When working with families of curves, the Eisenbud-Harris construction is forced to omit the locus of crude limit linear series. 
Brian Osserman recently introduced a new construction $G_d^r(X)$ parametrizing a modified notion of limit linear series \cite{osserman20061}. Osserman's construction is restricted to the case of curves with two components, but it allows the inclusion of crude limit linear series even when working with families of curves. 

When $X$ is a general reducible curve, in the Eisenbud-Harris construction, the set of refined limit linear series is always dense in $G_d^{r, \EH}(X)$. Osserman asked in \cite[Question 7.2]{osserman20061} whether the same is true in $G_d^r(X),$ and if not, which types of crude limit linear series contain an open subset of $G_d^r(X).$ The purpose of this paper is to answer this question. 

Accordingly, we suppose throughout that $X$ is a reducible curve with two smooth components $Y$ and $Z$ that intersect at the node $p.$ We recall the Eisenbud-Harris definition in this case. For any pair of linear series $V^Y$ and $V^Z$ on $Y$ and $Z$ respectively, both of dimension $r$ and degree $d,$ let $a^Y := \{a_i^Y\}_{i=0}^r$ and $a^Z := \{a_i^Z\}_{i=0}^r$ be the vanishing sequences of $V^Y$ and $V^Z$ at $p.$ By the definition of vanishing sequences, we have that
\begin{itemize}
\item[(A1)] $0 \le a_0^Y < a_1^Y < \cdots < a_r^Y \le d,$
\item[(A2)] $0 \le a_0^Z < a_1^Z < \cdots < a_r^Z \le d.$
\end{itemize}
The pair $(V^Y, V^Z)$ is an {\it Eisenbud-Harris (EH) limit linear series} of dimension $r$ and degree $d$ on $X$ if
\begin{itemize}
\item[(A3)] $a_i^Y + a_{r-i}^Z \ge d,$ for  $i = 0, 1, \dots, r.$
\end{itemize}
If equality holds in (A3) for all $i = 0,1,\dots,r,$ then $(V^Y, V^Z)$ is a {\it refined} EH limit linear series; otherwise, $(V^Y, V^Z)$ is a {\it crude} EH limit linear series.  

For any integer sequences $(a^Y, a^Z)$ satisfying (A1)--(A3), we denote by $G_d^r(X; a^Y, a^Z)$ the space consisting of EH limit linear series with vanishing sequences $(a^Y, a^Z)$. This gives a stratification of $G_d^{r, \EH}(X).$

We use the term {\it limit linear series} (without EH) to refer to the modified notion introduced by Osserman. For a fixed reducible curve $X,$ there is a surjective map $\pi$ from $G_d^r(X)$ to $G_d^{r, \EH}(X),$ and refined and crude limit linear series are defined as the preimage under $\pi$ of refined and crude EH limit linear series. Furthermore, $\pi$ is an isomorphism on refined limit linear series, but not on crude limit linear series, where it typically has positive-dimensional fibers. 

We need the following definition to state our main theorem. For convenience, we write $j \in a^Y$ if $j = a_i^Y$ for some $i,$ and write $j \in a^Z$ if $j = a_{i'}^Z$ for some $i'.$
\begin{defn}
Suppose $(a^Y, a^Z)$ satisfy (A1)--(A3). For any $i \in \{0, \dots, r\},$ we say $(a^Y, a^Z)$ is {\it connected at $i$} if there exists $j \in \{0, \dots, d\}$ such that $d -a_{r-i}^Z \le j \le a_i^Y,$ any integer between $j$ and $a_i^Y$ (inclusive) is in $a^Y$ and any integer between $d-j$ and $a_{r-i}^Z$ (inclusive) is in $a^Z.$ Sometimes we say $(a^Y, a^Z)$ is {\it connected at $i$ via $j$} if we want to indicate which $j$ we use.

$(a^Y, a^Z)$ is {\it connected} if $(a^Y, a^Z)$ is connected at every $i \in \{0, \dots, r\}.$
\end{defn}

We now state our main theorem.
\begin{thm}\label{main}
Suppose $X$ is a general reducible curve.
Then $\pi^{-1}(G_d^r(X; a^Y, a^Z))$ contains an open subset of $G_d^r(X)$ if and only if $G_d^r(X; a^Y, a^Z)$ is nonempty and $(a^Y, a^Z)$ is connected.
\end{thm}

Eisenbud and Harris showed that the dimension of each component of $G_d^{r, \EH}(X)$ is at least the ``Brill-Noether number'' $$\rho := (r+1)(d-r) - rg,$$ with equality if $X$ is a general reducible curve. In \cite{osserman20062}, Osserman was able to show this is also true for the space $G_d^r(X)$. Thus, when $X$ is general, $\pi^{-1}(G_d^r(X; a^Y, a^Z))$ contains an open subset of $G_d^r(X)$ if and only if it has dimension $\rho.$

In fact, without any generality hypothesis on $X,$ for each EH limit linear series $(V^Y, V^Z)$  with vanishing sequences $(a^Y, a^Z),$ Osserman gave a stratification of the fiber of $\pi$ over $(V^Y, V^Z)$, and also a purely combinatorial formula for the dimension of each stratum. Because of the complexity of these formulas, he did not analyze them directly. Instead, he gave an indirect argument to bound the dimension of fibers of $\pi$ in order to conclude that $\dim(G_d^r(X)) = \rho$ when $X$ is general.

In this paper, we use a combinatorial analysis of Osserman's formula to prove our main theorem by computing when we can have $\dim(\pi^{-1}(G_d^r(X; a^Y, a^Z))) = \rho.$

\section{Preliminaries}
In this section, we review Osserman's related work in detail, and state our problem in a combinatorial setting.
Recall $X$ is a reducible curve with two components $Y$ and $Z$ that intersect at the node $p.$ Suppose $(a^Y, a^Z)$ is the vanishing sequences of an EH limit linear series $(V^Y, V^Z)$ at $p.$ Thus, $(a^Y, a^Z)$ satisfies (A1)--(A3).

For $j = 0, \dots, d,$ let 
\begin{equation}
\tag{A2B}b_j^Y := \# \{ i : a_i^Y \ge j \}, \mbox{  and  } b_j^Z := \# \{ i : a_i^Z \ge d - j \}.
\end{equation}
We state without proof the following trivial lemma. 

\begin{lem}\label{bprop}
Suppose we have integers $a^Y = \{a_i^Y\}_{i=0}^r$ and $a^Z = \{a_i^Z\}_{i=0}^r$ satisfying (A1) and (A2), and integers $b^Y = \{ b_j^Y \}_{i=0}^d$ and $b^Z = \{ b_j^Z \}_{i=0}^d$ are given by (A2B). Then we have the following:
\begin{equation}\label{bYcmp}
\tag{B1} b_j^Y - b_{j+1}^Y = 
\begin{cases}1, & \mbox{if $j \in a^Y$}, \\ 0, & \mbox{if $j \not\in a^Y$,}\end{cases} 
\mbox{\ \ \    $\forall j = 0, 1, \dots, d-1.$}
\end{equation}
\begin{equation}\label{bZcmp}
\tag{B2} b_{j}^Z - b_{j-1}^Z = 
\begin{cases}1, & \mbox{if $d-j \in  a^Z$}, \\ 0, & \mbox{if $d-j \not\in a^Z$,} \end{cases} 
\mbox{\ \ \    $\forall j = 1, 2 \dots, d.$}
\end{equation}
\end{lem}

For each point of $G_d^r(X)$ in the fiber of $(V^Y, V^Z),$ Osserman associates integers 
$\{ \beta_j^Y\}_{j=0}^{d},$ $\{\beta_j^Z\}_{j=0}^{d},$ and $\{\epsilon_j\}_{j=1}^{d-1}$ satisfying for any $j = 1, \dots, d-1$
\begin{itemize}
\item[(C1)] $\epsilon_j = 0$ or $1,$ and $\epsilon_j = 0$ unless $j \in a^Y$ and $d-j \in a^Z.$
\item[(C2)] $b_{j+1}^Y \ge \beta_j^Y - \epsilon_j \ge \beta_{j+1}^Y.$
\item[(C3)] $b_{j-1}^Z \ge \beta_j^Z - \epsilon_j \ge \beta_{j-1}^Z.$
\item[(C4)] $\beta_j^Y + \beta_j^Z - \epsilon_j \ge r+1.$
\item[(C5)] $1 + r \ge \beta_{j+1}^Y + \beta_j^Z$ and $1 + r \ge \beta_j^Y + \beta_{j-1}^Z.$
\item[(C6)] $\beta_0^Y = \beta_d^Z = r+1,$ and $\beta_d^Y = \begin{cases}1, & a_0^Z = 0 \\ 0, & a_0^Z > 0\end{cases}$, $\beta_0^Z = \begin{cases}1, & a_0^Y = 0 \\ 0, & a_0^Y > 0\end{cases}$.
\end{itemize}
This gives a stratification of the fiber of $(V^Y, V^Z).$ Given $\beta^Y = \{ \beta_j^Y\}_{j=0}^{d},$ $\beta^Z = \{\beta_j^Z\}_{j=0}^{d},$ and $\epsilon = \{\epsilon_j\}_{j=1}^{d-1}$ satisfying the above conditions, Osserman \cite[Theorem 4.2, Lemma 5.4]{osserman20062} proves that the corresponding stratum is nonempty and has dimension given by
\begin{align*}
\D(\beta^Y, \beta^Z, \epsilon) := \sum_{j=1}^{d-1} &(\beta_j^Y - \beta_{j+1}^Y) (b_{j+1}^Y - \beta_j^Y + \epsilon_j) + (\beta_j^Z - \beta_{j-1}^Z)(b_{j-1}^Z - \beta_j^Z + \epsilon_j)  \\
& +  (r+1 - \beta_{j+1}^Y - \beta_{j-1}^Z)(\beta_j^Y + \beta_j^Z - \epsilon_j - r - 1).
\end{align*}

\begin{rem}
The dimension formula is given in (4.8) in \cite{osserman20062}. The notation there is different from what we use. The following table shows the correspondence between our notation and Osserman's in \cite[Theorem 4.2]{osserman20062}. \\
\begin{tabular}{|c|c|c|c|c|c|}
\hline
Indices & & & $1 \le j \le d$ & $0 \le j \le d-1$ & $1 \le j \le d-1$ \\
\hline
Our notation & $d$ & $r$  & $\beta_{j}^Y$ & $\beta_j^Z$ & $\epsilon_j$  \\
\hline
Osserman's notation & $n-1$ & $r-1$  & $d_{V_{1,j+1}}$ & $d_{V_{n,j+1}}$ & $d_{z_{j+1}}$  \\
\hline
\end{tabular}

In \cite[Theorem 4.2]{osserman20062}, Osserman gives restrictions (4.2)--(4.7) on $d_{V_{1,j+1}}$, $d_{V_{n,j+1}}$, and $d_{z_{i+1}}$, as well as conventions $d_{V_{1,n}} := d(V_{1,n})$  and $d_{V_{n,1}} := d(V_{n,1}),$ then concludes the dimension formula (4.8). We will verify that these conditions are equivalent to our conditions (C1)--(C6), and also that the dimension formulas are the same.

By Lemma 5.4 in \cite{osserman20062}, we have that 
\begin{alist}
\itm $b_j^Y = d(\overline{V}_{1,j+1})$ and $b_j^Z = d(\overline{V}_{n,j+1}),$ for $j = 0, 1 \dots, d.$ 
\itm $d(V_{1,n}) = \# \{i: a_i^Z \le 0 \}$ and $d(V_{n,1}) = \# \{i: a_i^Y \le 0 \}.$ 
\itm $d(\overline{Z}_{j+1}) = 1,$ for $j = 1, \dots, d-1.$
\end{alist}

First, it is clear that the conditions (C2)--(C5) are exactly the conditions (4.3)--(4.7) in \cite{osserman20062}. Second, by b), one sees that the second half of (C6) is equivalent to the convention $d_{V_{1,n}} := d(V_{1,n})$  and $d_{V_{n,1}} := d(V_{n,1}).$ Strictly speaking, Osserman does not define $\beta_0^Y$ and $\beta_d^Z,$ and they do not appear in the dimension formula, but we will find the convention $\beta_0^Y = \beta_d^Z = r+1$ convenient. Therefore, we have condition (C6).

Combining c) and \cite[Lemma 4.4/(i)]{osserman20062}, we have that
\begin{equation}\label{le1}
d(\overline{Z}_{1,j+1} \cap \overline{Z}_{n,j+1}) \le 1
\end{equation} 
and the equality holds if and only if $d(\overline{Z}_{1,j+1}) = d(\overline{Z}_{n,j+1}) = 1.$
Hence, the condition (4.2) in \cite{osserman20062} is equivalent to $\epsilon_j = 0$ or $1$ and $\epsilon_j =0$ unless $d(\overline{Z}_{1,j+1}) = d(\overline{Z}_{n,j+1}) = 1.$
However, by \cite[Lemma 3.4]{osserman20062}, $d(\overline{Z}_{1,j+1}) = d(\overline{V}_{1,j+1}) - d(\overline{V}_{1,j+2}) = b_j^Y - b_{j+1}^Y.$ Thus, by \eqref{bYcmp}, $d(\overline{Z}_{1,j+1}) =1$ if and only if $j \in a^Y.$ Similarly, $d(\overline{Z}_{n,j+1}) =1$ if and only if $d-j \in a^Z.$ Hence, condition (4.2) in \cite{osserman20062} is equivalent to (C1).

Finally, one notices that the only difference between the summands of $\D(\beta^Y, \beta^Z, \epsilon)$ and the dimension formula (4.8) in \cite{osserman20062} is that the latter has an additional term $d_{z_{j+1}}(d(\overline{Z}_{1,j+1} \cap \overline{Z}_{n,j+1}) - d_{z_{j+1}}).$ However, by condition (4.2) in \cite{osserman20062} and \eqref{le1},  in our situation this term is always zero.
\end{rem}

Via an indirect argument, Osserman also proves
\begin{thm}\cite[Corollary 5.5]{osserman20062}\label{main1}
For any $\beta^Y = \{ \beta_j^Y\}_{j=0}^{d},$ $\beta^Z = \{\beta_j^Z\}_{j=0}^{d},$ and $\epsilon = \{\epsilon_j\}_{j=1}^{d-1}$ satisfying conditions (C1)--(C6), we have that
\begin{equation}\label{mainineq}
\D(\beta^Y, \beta^Z, \epsilon) \le \sum_{i=0}^r (a_i^Y + a_{r-i}^Z - d).
\end{equation}
Hence, the dimension of the space of the crude limit series corresponding to a given EH limit linear series $(V^Y, V^Z)$ is bounded above by $\sum_{i=0}^r (a_i^Y + a_{r-i}^Z - d).$
\end{thm}

It is clear that the stratification and the dimension formula only depend on $(a^Y, a^Z)$, and are otherwise independent of the choice of $(V^Y, V^Z).$ Therefore, given any integer sequences $(a^Y, a^Z)$ satisfying (A1)--(A3) with nonempty $\pi^{-1}(G_d^r(X; a^Y, a^Z)),$ we have that
$$\dim(\pi^{-1}(G_d^r(X; a^Y, a^Z)) = \dim(\pi^{-1}(V^Y, V^Z)) + \dim(G_d^r(X; a^Y, a^Z)),$$
where $(V^Y, V^Z)$ is any EH limit linear series with vanishing sequences $(a^Y, a^Z).$ 

For $X$ general, Eisenbud and Harris \cite[Theorem 4.5]{eisenbudharris1983} showed that $$\dim(G_d^r(X; a^Y, a^Z)) = \rho - \sum_{i=0}^r (a_i^Y + a_{r-i}^Z - d).$$ 

We thus see that Theorem \ref{main} follows if we prove that for given $(V^Y, V^Z)$ with vanishing sequences $(a^Y, a^Z),$ we have 
$\dim(\pi^{-1}(V^Y, V^Z)) = \sum_{i=0}^r (a_i^Y + a_{r-i}^Z - d)$ if and only if $(a^Y, a^Z)$ is connected. However, $\dim(\pi^{-1}(V^Y, V^Z))$ is the maximum of $\D(\beta^Y, \beta^Z, \epsilon)$ over all possible $(\beta^Y, \beta^Z, \epsilon)$ satisfying (C1)--(C6). Hence, Theorem \ref{main} follows from the following theorem.

\begin{thm}\label{main2}
$(a^Y, a^Z)$ is connected if and only if there exist $\beta^Y = \{ \beta_j^Y\}_{j=0}^{d},$ $\beta^Z = \{\beta_j^Z\}_{j=0}^{d},$ and $\epsilon = \{\epsilon_j\}_{j=1}^{d-1}$ satisfying conditions (C1)--(C6), such that the equality in \eqref{mainineq} holds.
\end{thm}

Our goal for the rest of the paper is to prove Theorem \ref{main2}. In fact, we will also reprove Theorem \ref{main1} with a direct combinatorial analysis. In the process, we will establish a necessary and sufficient condition for the equality in \eqref{mainineq} to hold, which we use to prove Theorem \ref{main2}. 

\section{Synchronization of indices}
One difficulty in comparing the two sides of \eqref{mainineq} is that they are indexed differently. The left side $\D(\beta^Y, \beta^Z, \epsilon)$ is indexed by $j$ from $1$ to $d-1,$ and the right side is indexed by $i$ from $0$ to $r.$ In this section, we will discuss a relation between $i$'s and $j$'s, and use that to rewrite both side of \eqref{mainineq} so that they have the same indexing. For simplicity, we start with a definition.

Given any integer sequences $a^Y = \{a_i^Y\}_{i=0}^r$ and $a^Z = \{a_i^Z\}_{i=0}^r$ satisfying (A1)--(A3), integer sequences $b^Y = \{ b_j^Y \}_{i=0}^d$ and $b^Z = \{ b_j^Z \}_{i=0}^d$ are determined by (A2B). Thus, the conditions (C1)--(C6) are determined by $a^Y$ and $a^Z.$ If a set of integers $\beta^Y = \{ \beta_j^Y\}_{j=0}^{d},$ $\beta^Z = \{\beta_j^Z\}_{j=0}^{d},$ and $\epsilon = \{\epsilon_j\}_{j=1}^{d-1}$ satisfy conditions (C1)--(C6) associated to the given $a^Y$ and $a^Z,$ we call  $(\beta^Y, \beta^Z, \epsilon)$ {\it admissible with respect to $(a^Y, a^Z)$}. (We often omit ``with respect to $(a^Y, a^Z)$'' if there is no possibility of confusion.) 

From now on, we will assume that $(a^Y, a^Z)$ is a fixed pair of integer sequences satisfying (A1)--(A3), and $(b^Y, b^Z)$ are determined by (A2B). There are a few basic properties of an admissible $(\beta^Y, \beta^Z, \epsilon)$ that we often need to use in our proofs. We summarize them in the following lemma.

\begin{lem}Suppose $(\beta^Y, \beta^Z, \epsilon)$ is admissible. Then we have the following: 
\begin{equation}\label{bbetacmp}
b_j^Y \ge \beta_j^Y, b_j^Z \ge \beta_j^Z, \mbox{\ \ \    $\forall j = 0,1, \dots, d.$}
\end{equation}
\begin{equation}\label{betaY}
r+1= \beta_0^Y \ge \beta_1^Y \ge \cdots \ge \beta_d^Y \ge 0.
\end{equation}
\begin{equation}\label{betaZ}
0 \le \beta_0^Z \le \beta_1^Z \le \cdots \le \beta_d^Z = r+1.
\end{equation}
\begin{equation}\label{betaZneq0}
\mbox{If $\beta_0^Z \neq 0,$ then $\beta_0^Z = 1,$ $a_0^Y = 0$ and $a_r^Z = d.$}
\end{equation}
\begin{equation}\label{betaYneq0}
\mbox{If $\beta_d^Y \neq 0,$ then $\beta_d^Y = 1,$ $a_0^Z = 0$ and $a_r^Y = d.$}
\end{equation}
\begin{equation}\label{zero1}
b_{1}^Y + \beta_{0}^Z - r-1 = 0.
\end{equation}
\begin{equation}\label{zero2}
b_{d-1}^Z + \beta_{d}^Y - r-1= 0.
\end{equation}
\end{lem}

\begin{proof}
By  (C2), it is clear that $b_j^Y \ge \beta_j^Y$ for any $j = 2, \dots, d.$ If $j = 0,$ by (A2B) and (C6), $b_j^Y = r+1 = \beta_j^Y.$ If $j = 1,$ by (C2), \eqref{bYcmp} and (C1), we have that $\beta_j^Y \le b_{j+1}^Y + \epsilon_j \le b_j^Y.$ Therefore, $b_j^Y \ge \beta_j^Y$ for any $j = 0, 1, \dots, d.$ Similarly, we can show that $b_j^Z \ge \beta_j^Z$ for any $j = 0, 1, \dots, d.$ Thus, \eqref{bbetacmp} holds.

$\beta_0^Y \ge \beta_1^Y$ is the only part of \eqref{betaY} that does not directly follow from conditions (C1)--(C6). However, by \eqref{bbetacmp}, we have $\beta_1^Y \le b_1^Y \le r+1 = \beta_0^Y.$ Thus, \eqref{betaY} holds. We can similarly show that \eqref{betaZ} is true.

Both \eqref{betaZneq0} and \eqref{betaYneq0} follow from (C6) and (A1)--(A3).
One can verify \eqref{zero1} and \eqref{zero2} by considering whether $0 \in a^Y$ and whether $0 \in a^Z,$ respectively.
\end{proof}

The following lemma is the first step of our rewriting of the right hand side of \eqref{mainineq}.
\begin{lem}
Suppose $(\beta^Y, \beta^Z, \epsilon)$ is admissible. Then
\begin{equation}\label{rewriteR}
\sum_{i=0}^r (a_i^Y + a_{r-i}^Z - d) = \sum_{i=\beta_0^Z}^{r} (a_i^Y + a_{r-i}^Z - d) = \sum_{i=\beta_0^Z}^{r-\beta_d^Y} (a_i^Y + a_{r-i}^Z - d)
\end{equation}
\end{lem}
\begin{proof}
The first equality in \eqref{rewriteR} holds if $\beta_0^Z = 0.$ Suppose $\beta_0^Z \neq 0.$ 
Then by \eqref{betaZneq0}, $a_0^Y = 0$ and $a_r^Z = d,$ so $a_0^Y + a_r^Z - d = 0.$ 
Hence, the first equality in \eqref{rewriteR} still holds. We can similarly show that the second equality in \eqref{rewriteR} always holds.
\end{proof}

Because of the above lemma, we will focus on $i$ between $\beta_{0}^Z$ and $r-\beta_d^Y.$
Given an admissible $(\beta^Y, \beta^Z, \epsilon)$, for any $i$ with $\beta_{0}^Z \le i \le r-\beta_d^Y,$ by (C6), we have
$$r+1 - \beta_0^Y \le \beta_{0}^Z \le i \le r-\beta_d^Y \le \beta_d^Z-1.$$
Thus, by \eqref{betaY} and \eqref{betaZ}, there exist unique $j_1 \in \{0, \dots, d-1\}$ and $j_2 \in \{1, \dots, d\}$ such that 
\begin{equation}\label{boundY}
r+1 - \beta_{j_1}^Y \le i \le r - \beta_{j_1+1}^Y,
\end{equation}
\begin{equation}\label{boundZ}
\beta_{j_2-1}^Z \le i \le \beta_{j_2}^Z - 1.
\end{equation}

\begin{defn}
Suppose $(\beta^Y, \beta^Z, \epsilon)$ is admissible. Writing $\beta = (\beta^Y, \beta^Z),$ we define
$$I_\beta = I_{(\beta^Y, \beta^Z)} := \{i: \beta_{0}^Z \le i \le r-\beta_d^Y\},$$
and then define a map
$$\psi_\beta = \psi_{(\beta^Y, \beta^Z)}: I_\beta \to \{0, \dots, d-1\} \times \{1, \dots, d\},$$
by mapping each $i \in I_\beta$ to the pair $(j_1, j_2)$ determined by \eqref{boundY} and \eqref{boundZ}.
\end{defn}

\begin{lem}\label{psibeta}
Suppose $(\beta^Y, \beta^Z, \epsilon)$ is admissible. For any $i \in I_\beta,$ $\psi_\beta(i)$ is either $(j,j),$ for some $j \in \{1, \dots, d-1\},$ or $(j-1,j),$ for some $j \in \{1, \dots, d\}.$ Furthermore,
\begin{equation*}
\psi_\beta(i) = \begin{cases}(j, j) \mbox{ for some } j \in \{1, \dots, d-1\}, &  \mbox{ if and only if $r+1 - \beta_j^Y \le i \le \beta_j^Z -1$}; \\
(j-1,j) \mbox{ for some } j \in \{1, \dots, d\}, & \mbox{ if and only if $\beta_{j-1}^Z \le i \le r - \beta_j^Y.$}
\end{cases}
\end{equation*}
\end{lem}

\begin{proof}
Let $\psi_\beta(i) = (j_1, j_2).$ We first show that $j_2 - 1 \le j_1 \le j_2.$ Since \eqref{boundY} and \eqref{boundZ} are both satisfied,  we must have
$$r+1 - \beta_{j_1}^Y \le \beta_{j_2}^Z - 1 \mbox{ and } \beta_{j_2-1}^Z \le r - \beta_{j_1+1}^Y.$$

Thus, if $j_1 \neq 0,$ then by using (C5) and \eqref{betaZ}, 
$$\beta_{j_1-1}^Z \le r+1 - \beta_{j_1}^Y < \beta_{j_2}^Z \Rightarrow j_1 -1 < j_2 \Rightarrow j_1 \le j_2.$$ 
If $j_1 = 0,$ it is clear that we still have $j_1 \le j_2.$

If $j_1 \neq d-1,$ then by using (C4) and \eqref{betaZ},
$$\beta_{j_2-1}^Z < r+1 - \beta_{j_1+1}^Y \le \beta_{j_1+1}^Z \Rightarrow j_2-1 < j_1+1 \Rightarrow j_2-1 \le j_1.$$
If $j_1 = d-1,$ we also have $j_2 - 1 \le j_1.$

Hence, we must have that $j_1 =j_2$ or $j_1 = j_2-1.$ Therefore, $\psi_\beta(i) = (j,j)$ or $(j-1,j)$ for some $j.$

By the definition of $\psi_\beta,$ we have that $\psi_\beta(i) = (j, j)$ if and only if \eqref{boundY} and \eqref{boundZ} are satisfied for $j_1 = j_2 = j,$ which is equivalent to
$$r+1 - \beta_{j}^Y = \max(r+1 - \beta_{j}^Y, \beta_{j-1}^Z) \le i \le \min(r - \beta_{j+1}^Y, \beta_j^Z-1) = \beta_j^Z-1.$$
Also, since $j_1 \in \{0, \dots, d-1\}$ and $j_2 \in \{1, \dots, d\},$ we must have that $j \in \{1, \dots, d-1\}.$

Similarly, we can show that $\psi_\beta(i) = (j-1,j)$ for some $j \in \{1, \dots, d\}$ if and only if $\beta_{j-1}^Z \le i \le r - \beta_j^Y.$
\end{proof}

\begin{cor}\label{preimg}
Suppose $(\beta^Y, \beta^Z, \epsilon)$ is admissible.
For any $j \in \{1, \dots, d-1\},$ 
$$\# \psi_\beta^{-1}(j,j) = \# \{ i \in I_\beta : \psi_\beta(i) = (j, j) \} = \beta_j^Y + \beta_j^Z - (r+1).$$
Thus, $\psi_\beta^{-1}(j,j)$ is nonempty if and only if $\beta_j^Y + \beta_j^Z \neq r+1.$

For any $j \in \{1, \dots, d\},$
$$\# \psi_\beta^{-1}(j-1,j) =\# \{ i \in I_\beta : \psi_\beta(i) = (j-1, j) \} = r+1 - \beta_j^Y - \beta_{j-1}^Z.$$
Thus, $\psi_\beta^{-1}(j-1,j)$ is nonempty under $\psi_\beta$ if and only if $\beta_j^Y + \beta_{j-1}^Z \neq r+1.$

\end{cor}

\begin{proof}
We only need to check that for any $j \in \{1, \dots, d-1\},$  any $i$ satisfying $r+1 - \beta_j^Y \le i \le \beta_j^Z -1$ is in $I_\beta,$ and for any $j \in \{1, \dots, d\},$ any $i$ satisfying $\beta_{j-1}^Z \le i \le r - \beta_j^Y$ is in $I_\beta.$ Both can be verified easily.
\end{proof}

\begin{defn}
Suppose $(\beta^Y, \beta^Z, \epsilon)$ is admissible. Define
$$J_\beta = J_{(\beta^Y, \beta^Z)} := \{ j: 1 \le j \le d-1 \ : \  \psi_\beta^{-1}(j,j) \neq \emptyset \},$$
$$J'_\beta = J'_{(\beta^Y, \beta^Z)} := \{ j: 1 \le j \le d  \ : \ \psi_\beta^{-1}(j-1,j) \neq \emptyset\}.$$
\end{defn}

Note that the sets $J_\beta$ and $J'_\beta$ are not necessarily disjoint. Now we state the main lemma of this section.
\begin{lem}\label{ineqpsi}
Suppose $(\beta^Y, \beta^Z, \epsilon)$ is admissible. Then we have the following two identities:
\begin{multline}
\sum_{i=0}^r (a_i^Y + a_{r-i}^Z - d) = \sum_{j \in J_\beta} \sum_{i \in \psi_\beta^{-1}(j,j)} (a_i^Y + a_{r-i}^Z - d)  + \sum_{j \in J'_\beta} \sum_{i \in \psi_\beta^{-1}(j-1,j)} (a_i^Y + a_{r-i}^Z - d), \label{Rwpsi}
\end{multline}
\begin{align}
\D(\beta^Y, \beta^Z, \epsilon) =&  \sum_{j \in J_\beta} \sum_{i \in \psi_\beta^{-1}(j,j)} (b_{j+1}^Y + b_{j-1}^Z + \epsilon_j - r -1) + \label{Lwpsi}\\
&  \sum_{j \in J'_\beta} \sum_{i \in \psi_\beta^{-1}(j-1,j)} \left( (b_{j}^Y + \beta_{j-1}^Z - r-1)+ (b_{j-1}^Z + \beta_{j}^Y - r-1)\right). \nonumber
\end{align}
\end{lem}

\begin{proof}
\eqref{Rwpsi} follows immediately from \eqref{rewriteR} and Lemma \ref{psibeta}. Thus, we only need to prove \eqref{Lwpsi}.
We first rearrange the summand in the definition of $\D(\beta^Y, \beta^Z, \epsilon)$ and get
\begin{multline*}
\D(\beta^Y, \beta^Z, \epsilon) = \sum_{j=1}^{d-1} (\beta_j^Y + \beta_j^Z - r-1) (b_{j+1}^Y + b_{j-1}^Z + \epsilon_j - r -1) + \\
  \hspace{5mm} \sum_{j=1}^{d-1} (r+1 - \beta_{j+1}^Y - \beta_{j}^Z) (b_{j+1}^Y + \beta_{j}^Z - r-1)+  \sum_{j=1}^{d-1} (r+1 - \beta_j^Y - \beta_{j-1}^Z) (b_{j-1}^Z + \beta_{j}^Y - r-1).
\end{multline*}
Then by using \eqref{zero1} and \eqref{zero2}, we rewrite the second and third sums in the above formula 
\begin{multline*}
 \sum_{j=1}^{d-1} (r+1 - \beta_{j+1}^Y - \beta_{j}^Z) (b_{j+1}^Y + \beta_{j}^Z - r-1)+  \sum_{j=1}^{d-1}(r+1 - \beta_j^Y - \beta_{j-1}^Z) (b_{j-1}^Z + \beta_{j}^Y - r-1) \\
= \sum_{j=0}^{d-1} (r+1 - \beta_{j+1}^Y - \beta_{j}^Z) (b_{j+1}^Y + \beta_{j}^Z - r-1)+ \sum_{j=1}^d (r+1- \beta_j^Y - \beta_{j-1}^Z) (b_{j-1}^Z + \beta_{j}^Y - r-1) \\
\shoveleft{=  \sum_{j=1}^d (r+1 - \beta_j^Y - \beta_{j-1}^Z) \left( (b_{j}^Y + \beta_{j-1}^Z - r-1)+ (b_{j-1}^Z + \beta_{j}^Y - r-1)\right). }   
\end{multline*}
Therefore,
\begin{align*}
\D(\beta^Y, \beta^Z, \epsilon) &= \sum_{j=1}^{d-1} (\beta_j^Y + \beta_j^Z - r-1) (b_{j+1}^Y + b_{j-1}^Z + \epsilon_j - r -1) + \\
&   \sum_{j=1}^d (r+1 - \beta_j^Y - \beta_{j-1}^Z) \left( (b_{j}^Y + \beta_{j-1}^Z - r-1)+ (b_{j-1}^Z + \beta_{j}^Y - r-1)\right).
\end{align*}
However, by Corollary \ref{preimg}, we know that $\beta_j^Y + \beta_j^Z - (r+1)=0$ for any $j \in \{1, \dots,d-1\} \setminus J_\beta$, and $(r+1) - \beta_j^Y - \beta_{j-1}^Z=0$ for any $j \in \{1, \dots, d\} \setminus J'_\beta.$ Thus,
\begin{multline*}
\D(\beta^Y, \beta^Z, \epsilon) = \sum_{j \in J_\beta} (\beta_j^Y + \beta_j^Z - r-1) (b_{j+1}^Y + b_{j-1}^Z + \epsilon_j - r -1) + \\
 \sum_{j \in J'_\beta} (r+1 - \beta_j^Y - \beta_{j-1}^Z) \left( (b_{j}^Y + \beta_{j-1}^Z - r-1)+ (b_{j-1}^Z + \beta_{j}^Y - r-1)\right).
\end{multline*}
Applying Corollary \ref{preimg} again, we obtain \eqref{Lwpsi}.
\end{proof}

\section{Comparison of terms}
Lemma \ref{ineqpsi} synchronized the indices on both sides of \eqref{mainineq}. In this section, we will show each term in \eqref{Lwpsi} is less than or equal to the corresponding term in \eqref{Rwpsi}, then conclude Proposition \ref{mainstr}, which is a stronger version of Theorem \ref{main1} including a necessary and sufficient condition for the equality in \eqref{mainineq} to hold.

\begin{lem}\label{j1j2prop}
Suppose $(\beta^Y, \beta^Z, \epsilon)$ is admissible. Let $i \in I_\beta,$ and $\psi_\beta(i) = (j_1, j_2).$ Then we have the following:
\begin{ilist}
\itm $a_{i}^Y \ge j_1,$ and $a_{i}^Y - j_1 \ge b_{j_1+1}^Y - (r-i),$ where the second equality holds if and only if either $a_{i}^Y = j_1,$ or $a_{i}^Y > j_1$ and any integer between $j_1+1$ and $a_{i}^Y$ (inclusive) is in $a^Y$.
\itm $a_{r-i}^Z \ge d-j_2,$ and $a_{r-i}^Z - (d-j_2) \ge b_{j_2-1}^Z - i,$ where the second equality holds if and only if either $a_{r-i}^Z = d-j_2,$  or $a_{r-i}^Z > d-j_2$ and any integer between $d-j_2+1$ and $a_{r-i}^Z$ (inclusive) is in $a^Z$.
\end{ilist}
\end{lem}

\begin{proof}
By \eqref{boundY}, \eqref{bbetacmp} and (A2B), we have that
$$i \ge r+1-\beta_{j_1}^Y \ge r+1 - b_{j_1}^Y = \# \{i_0 \  : \  a_{i_0}^Y < j_1 \}.$$
Thus, $a_{i}^Y \ge j_1.$ If $a_{i}^Y = j_1,$ then $j_1 \in a^Y.$ Thus, by \eqref{bYcmp}, $b_{a_{i}^Y}^Y = b_{j_1}^Y = 1+b_{j_1+1}^Y.$ Thus, $a_{i}^Y - j_1 = 0 = 1+b_{j_1+1}^Y - b_{a_{i}^Y}^Y.$ If $a_{i}^Y > j_1,$ then $a_{i}^Y \ge j_1+1.$ Hence, 
\begin{eqnarray*}
a_{i}^Y - j_1 &=& 1 + a_{i}^Y - (j_1+1) \\
&\ge& 1 + \# \{i_0 : j_1+1 \le a_{i_0}^Y < a_{i}^Y \} \\
&=& 1 + \# \{i_0 :  a_{i_0}^Y \ge j_1+1\} - \#\{ i_0 : a_{i_0}^Y  \ge a_{i}^Y \} \\
&=&1+b_{j_1+1}^Y - b_{a_{i}^Y}^Y
\end{eqnarray*}
Note that by (A1), $b_{a_{i}^Y}^Y = (r+1) - i.$ Therefore, (i) is true.

Similarly, we can prove (ii).
\end{proof}

Now we can prove the desired inequalities between the terms in \eqref{Lwpsi} and \eqref{Rwpsi}.
\begin{lem}\label{lbd}
Suppose $(\beta^Y, \beta^Z, \epsilon)$ is admissible. Let $i \in I_\beta.$
\begin{ilist}
\itm If $\psi_\beta(i) = (j,j)$ for some $j \in \{1, \dots, d-1\},$ then 
\begin{equation}\label{suma-d1}
a_{i}^Y + a_{r-i}^Z - d \ge b_{j+1}^Y   + b_{j-1}^Z  + \epsilon_j- r - 1,
\end{equation}
where the equality holds if and only if $(a^Y, a^Z)$ is connected at $i$ via $j$ and $\epsilon_j =1.$

\itm If $\psi_\beta(i) = (j-1,j)$ for some $j \in \{1, \dots, d\},$ then
\begin{equation}\label{suma-d2}
a_{i}^Y + a_{r-i}^Z - d > (b_{j}^Y + \beta_{j-1}^Z - r-1)+ (b_{j-1}^Z + \beta_{j}^Y - r-1).
\end{equation}
\end{ilist}
\end{lem}

\begin{proof}
Suppose $\psi_\beta(i) = (j,j).$ Combining (i) and (ii) of Lemma \ref{j1j2prop} together, setting $j_1=j_2=j$ and using (C1),  we and
\begin{eqnarray*}
a_{i}^Y + a_{r-i}^Z - d &=& (a_{i}^Y - j_1)  +  (a_{r-i}^Z - (d - j_2)) \\
&\ge& b_{j_1+1}^Y - (r-i) + b_{j_2-1}^Z - i \\
&=& b_{j+1}^Y + b_{j-1}^Z - r \\
&\ge& b_{j+1}^Y + b_{j-1}^Z + \epsilon_j - r -1,
\end{eqnarray*}
where the equality holds if and only if any integer between $j$ and $a_{i}^Y$ (inclusive) is in $a^Y$, any integer between $d-j$ and $a_{r-i}^Z$ (inclusive) is in $a^Z,$ and $\epsilon_j =1.$ Note that Lemma \ref{j1j2prop} also gives that $d -a_{r-i}^Z \le j \le a_i^Y.$ Therefore, we conclude that equality holds in \eqref{suma-d1} if and only if $(a^Y, a^Z)$ is connected at $i$ via $j$ and $\epsilon_j =1.$

Suppose $\psi_\beta(i) = (j-1,j).$ Then $i \in \{ i_0 \in I_\beta : \psi_\beta(i_0) = (j-1, j) \}.$  By Corollary \ref{preimg}, $$(r+1) - \beta_j^Y - \beta_{j-1}^Z > 0.$$
Combining (i) and (ii) of Lemma \ref{j1j2prop} together, setting $(j_1,j_2) = (j-1,j)$ and using the above inequality, we get
\begin{eqnarray*}
a_{i}^Y + a_{r-i}^Z - d &=& (a_{i}^Y - j_1)  +  (a_{r-i}^Z - (d - j_2)) -1\\
&\ge& b_{j_1+1}^Y - (r-i) + b_{j_2-1}^Z - i -1\\
&=& b_{j}^Y + b_{j-1}^Z - r -1\\
&>& (b_{j}^Y + \beta_{j-1}^Z - r-1)+ (b_{j-1}^Z + \beta_{j}^Y - r-1).
\end{eqnarray*}
\end{proof}

The following proposition, which is the main result of the section, follows immediately from Lemma \ref{ineqpsi} and Lemma \ref{lbd}.
\begin{prop}\label{mainstr}
Suppose $(\beta^Y, \beta^Z, \epsilon)$ is admissible. 
Then
\begin{equation}\label{mainineq2}
\D(\beta^Y, \beta^Z, \epsilon) \le  \sum_{i=0}^r (a_i^Y + a_{r-i}^Z - d),
\end{equation}
where the equality holds if and only if both of the following two conditions hold:
\begin{ilist}
\itm $J'_\beta$ is the empty set, or equivalently (by Corollary \ref{preimg}), 
\begin{equation*}
\beta_j^Y + \beta_{j-1}^Z = r+1, \mbox{ for any $j = 1, \dots, d.$}
\end{equation*}
\itm For any $j \in J_\beta$ and any $i \in \psi_\beta^{-1}(j,j),$ we have that 
$(a^Y, a^Z)$ is connected at $i$ via $j,$
and $\epsilon_j =1.$
\end{ilist}
\end{prop}


Theorem \ref{main1} is an immediate corollary to Proposition \ref{mainstr}.

\begin{cor}\label{nece}
If there exists an admissible $(\beta^Y, \beta^Z, \epsilon)$ such that equality holds in \eqref{mainineq2}, then $(a^Y, a^Z)$ is connected.
\end{cor}

\begin{proof}
Suppose $(\beta^Y, \beta^Z, \epsilon)$ is admissible. We must have that conditions (i) and (ii) of Proposition \ref{mainstr} hold. $J'_\beta$ is empty, so for any  $i \in I_\beta,$ $\psi_\beta(i) = (j, j)$ for some $j \in \{1, \dots, d-1\}.$ Thus, $(a^Y, a^Z)$ is connected at any $i \in I_\beta.$ Recall $I_\beta = \{ \beta_0^Z \le i \le r - \beta_d^Y \}.$ However, by (C6), the only possible $i$'s not belonging to $I_\beta$ are $0$ and $d.$ If $0 \not\in I_\beta,$ then $\beta_0^Z = 1.$ Thus, by \eqref{betaZneq0}, we have that $a_0^Y = 0$ and $a_r^Z = d.$ Hence, $(a^Y, a^Z)$ is connected at $i=0$ via $j = 0.$ We can similarly show that if $r \not\in I_\beta,$ then $(a^Y, a^Z)$ is connected at $i=r$ via $j = d.$
\end{proof}

\section{Proof of the main theorems}
In the previous section, we showed that the connectedness of $(a^Y, a^Z)$ is a necessary condition for having an admissible $(\beta^Y, \beta^Z, \epsilon)$ which achieves equality in \eqref{mainineq}. In this section, we will show it is also a sufficient condition. 

\begin{prop}\label{suff}
If $(a^Y, a^Z)$ is connected, then there exists an admissible $(\beta^Y, \beta^Z, \epsilon)$ such that 
$$\D(\beta^Y, \beta^Z, \epsilon) = \sum_{i=0}^r (a_i^Y + a_{r-i}^Z - d).$$
\end{prop}

The idea of the proof is constructive, i.e., we will explicitly describe how to construct an admissible $(\beta^Y, \beta^Z, \epsilon)$ from a connected $(a^Y, a^Z)$ such that the equality holds. We will discuss properties of a connected $(a^Y, a^Z)$ and use those to make our construction.

\begin{lem}\label{cnnt1}
Suppose $(a^Y, a^Z)$ is connected. For any $0 \le i_1 \le i_2 \le r,$ if $(a^Y, a^Z)$ is connected at $i_1$ and $i_2$ via $j_1$ and $j_2$, respectively, with $j_1 \ge j_2,$ then $(a^Y, a^Z)$ is connected at $i$ via $j,$ for any $i: i_1 \le i \le i_2$ and any $j: j_1 \ge j \ge j_2.$
\end{lem}

\begin{proof}
Since $(a^Y, a^Z)$ is connected at $i_1$ and $i_2$ via $j_1$ and $j_2$, respectively,  we have that  $j_1 \le a_{i_1}^Y$, 
and any integer between $j_2$ and $a_{i_2}^Y$ is in $a^Y.$ However, because $j_2 \le j \le j_1$ and by (A1), we have that $a_{i_1}^Y \le a_i^Y \le a_{i_2}^Y,$ we conclude that
$$j_2 \le j \le j_1 \le a_{i_1}^Y \le a_i^Y \le a_{i_2}^Y,$$
so any integer between $j$ and $a_{i}^Y$ is in $a^Y.$
Similarly, we can show that $d -a_{r-i}^Z \le j$ and any integer between $d-j$ and $a_{r-i}^Z$ is in $a^Z.$ Hence, $(a^Y, a^Z)$ is connected at $i$ via $j.$ 
\end{proof}

\begin{lem}\label{cnnt2}
Suppose $(a^Y, a^Z)$ is connected. There exist $\nHat{J} \subseteq \{0, \dots, d\}$ and a surjection
$$\j: \{i: 0 \le i \le r\} \to \nHat{J}$$ satisfying
\begin{alist}
\itm For any $i: 0 \le i \le r,$ $(a^Y, a^Z)$ is connected at $i$ via $\j(i).$
\itm For any $0 \le i_1 < i_2 \le r,$ we have that $\j(i_1) \le \j(i_2).$
\end{alist}
\end{lem}

\begin{proof}
For any $i:0 \le i \le r,$ we define $\j(i)$ to be the greatest $j$ such that $(a^Y, a^Z)$ is connected at $i$ via $j.$ Let $\nHat{J} = \{ \j(i) \ : \ 0 \le i \le r\}.$ We check that $\nHat{J}$ and $\j: i \mapsto \j(i)$ satisfy the required conditions. First, it is clear that $\j: i \mapsto \j(i)$ is a surjection, and a) holds by the definition of $\j(i).$ Suppose b) does not hold. Then there exist $0 \le i_1 < i_2 \le r$ such that $\j(i_1) > \j(i_2).$ However, by Lemma \ref{cnnt1}, $(a^Y, a^Z)$ is connected at $i_2$ via $\j(i_1),$ which contradicts the definition of $\j(i_2).$
\end{proof}

\begin{lem}\label{cnnt3}
Suppose $(a^Y, a^Z)$ is connected. Let $\nHat{J}$ and $\j$ satisfy the conditions of Lemma \ref{cnnt2}.  Define $J = \nHat{J} \setminus \{0, d\}.$ Suppose the elements in $J$ are $j_1 < j_2 < \dots < j_s.$ Set $j_0 = 0$ and $j_{s+1}=d.$
For any $k: 0 \le k \le s+1,$ we define
\begin{equation}\label{Ikdefn}
I_k = \{i \ | \ \j(i) \ge j_k\}.
\end{equation} 
Then we have the following
\begin{equation}\label{j(i)Ik}
i \ge r+1 - |I_k| \Leftrightarrow \j(i) \ge j_k, \mbox{ for any $k = 0, \dots, s+1.$}
\end{equation}
\begin{equation}\label{bYIk}
b_{j_k}^Y \ge |I_k|, \mbox{ for any $k = 0, \dots, s+1.$}
\end{equation}
\begin{equation}\label{bZIk}
b_{j_k}^Z \ge r+1 - | I_{k+1}|, \mbox{ for any $k = 0, \dots, s.$}
\end{equation}
\begin{equation}\label{Ikdiff}
|I_k| - |I_{k+1}| \ge 1, \mbox{ for any $k = 1, \dots, s.$}
\end{equation}
\begin{equation}\label{0andr}
 I_{s+1}= \begin{cases}\{r\}, & a_0^Z = 0 \\ \emptyset, & a_0^Z > 0\end{cases}, 
\mbox{ and }
\{0, \dots, r\} \setminus I_1 = \begin{cases}\{0\}, & a_0^Y = 0 \\ \emptyset, & a_0^Y > 0\end{cases}.
\end{equation}
\end{lem}

\begin{proof}
Because of condition b) described in Lemma \ref{cnnt2}, the $i$'s in $I_k$ must be the largest $|I_k|$ $i$'s in $\{0, \dots, r\}.$ Thus, $i \in I_k$ if and only if $r+1 - |I_k| \le i \le r.$ Then \eqref{j(i)Ik} follows .

For any $i \in I_k,$ since $a_i^Y \ge \j(i) \ge j_k,$ $i \in \{i_0: a_{i_0}^Y \ge j_k\}.$ Thus, $I_k \subseteq \{i_0: a_{i_0}^Y \ge j_k\},$ so $b_{j_k}^Y \ge |I_k|.$

Let $k \in \{0, \dots, s\}.$ For any $i \not\in I_{k+1},$ we have that $\j(i) < j_{k+1}.$ Because $\j(i) \in \nHat{J} \subseteq \{j_0, j_1, \dots, j_s, j_{s+1}\},$ we know that $\j(i) \le j_k.$ Thus, $a_{r-i}^Z \ge d - \j(i) \ge d - j_k.$ Hence, $i \in \{ i_0 \ : \ a_{r-i_0}^Z \ge d-j_k\}.$ We conclude that $\{0,\dots, r\} \setminus I_{k+1} \subseteq \{ i_0 \ : \ a_{r-i_0}^Z \ge d-j_k\}.$ Therefore, $b_{j_k}^Z \ge r+1 - | I_{k+1}|.$

For any $k \in \{1, \dots, s\},$ since $j_k \in \nHat{J},$ by the definition of $\nHat{J},$ there exists $i$ such that $\j(i) = j_k.$ Thus, $I_k \setminus I_{k+1}$ is nonempty, and \eqref{Ikdiff} follows.

$I_{s+1} = \{i: \j(i) \ge d\} = \{i: \j(i) = d\}.$ For any $i$ with $\j(i) = d,$ we have that $a_i^Y \ge d,$ which implies $a_i^Y = d$ and $i = r.$ Thus, $I_{s+1} \subseteq \{r\}.$ However, one checks that there exist $i$ such that $\j(i) = d$ if and only if $a_0^Z = 0.$ 
Therefore, the first half of \eqref{0andr} holds.
We can similarly show the second half of \eqref{0andr} holds as well.
\end{proof}

\begin{proof}[Proof of Proposition \ref{suff}]

Suppose $(a^Y, a^Z)$ is connected. Let $\nHat{J}$ and $\j$ satisfy the conditions of Lemma \ref{cnnt2}. Define $J,$ $j_0,$ $j_1,$ $\dots,$ $j_{s+1},$ $I_0,$ $I_1,$ $\dots,$ $I_{s+1}$ as in Lemma \ref{cnnt3}. 

Note that for any $j$ with $1 \le j \le d,$ since $0 = j_0 < j_1 < \cdots < j_s < j_{s+1} = d,$ there exists a unique $k \in \{1, \dots, s+1\}$ such that $j_{k-1} < j \le j_k.$ We define $\beta^Y = \{ \beta_j^Y\}_{j=0}^{d},$ $\beta^Z = \{\beta_j^Z\}_{j=0}^{d},$ and $\epsilon = \{\epsilon_j\}_{j=1}^{d-1}$ as follows:
\begin{equation}\label{0andd}
\beta_0^Y = \beta_d^Z = r+1,
\end{equation}
\begin{equation}
\epsilon_j = \begin{cases}
1, & \mbox{ if $j \in J$,} \\ 0, & \mbox{ if $j \not\in J$,} 
\end{cases}
\end{equation}
and for any $j: 1 \le j \le d,$
\begin{equation}\label{betadefn}
\beta_j^Y = |I_k|, \mbox{ and } \beta_{j-1}^Z = r + 1 - |I_k|, \mbox{ if $j_{k-1} < j \le j_k.$}
\end{equation}

We now show that $(\beta^Y, \beta^Z, \epsilon)$ is admissible. For any $j \in J,$ there exists $i$ such that $\j(i) = j.$ Thus, by the definition of $(a^Y, a^Z)$ being connected at $i$ via $j,$ we have that $j \in a^Y$ and $d-j \in a^Z.$ Therefore, (C1) is satisfied. From \eqref{betadefn}, one sees that
\begin{equation}\label{r+1cond}
\beta_j^Y + \beta_{j-1}^Z = r+1, \mbox{ for any $j = 1, \dots, d.$}
\end{equation}
Therefore, (C5) is satisfied. 
By \eqref{0andr} and \eqref{0andd}, it is clear that (C6) holds. Below, we verify that (C2)--(C4) are true.

Let $j \in \{1, \dots, d-1\}.$  We first consider the case when $j_{k-1} < j < j_k$ for some $k \in \{1, \dots, s+1\}.$ Then by \eqref{bYcmp}, \eqref{bZcmp}, \eqref{bYIk} and \eqref{bZIk},
$$b_{j+1}^Y \ge b_{j_k}^Y \ge |I_k| = \beta_j^Y = \beta_j^Y - \epsilon_j = \beta_{j+1}^Y,$$
\begin{equation}\label{j_j-1Z1}
b_{j-1}^Z \ge b_{j_{k-1}}^Z \ge r+1 - |I_k| =  \beta_{j}^Z = \beta_{j}^Z - \epsilon_j = \beta_{j-1}^Z,
\end{equation}
$$\beta_j^Y + \beta_j^Z - \epsilon_j = \beta_{j+1}^Y + \beta_j^Z = r+1.$$
Hence, (C2)--(C4) are satisfied in this case.

Now suppose $j = j_k,$ for some $k \in \{1, \dots, s\}.$ (Note that $k \neq s+1$ because $j < d.$) Then by \eqref{bYcmp}, \eqref{bZcmp}, \eqref{bYIk}, \eqref{bZIk}, \eqref{Ikdiff}
$$b_{j+1}^Y = b_j^Y - 1 \ge |I_k| - 1 = \beta_j^Y - \epsilon_j \ge |I_{k+1}| = \beta_{j+1}^Y,$$
\begin{equation}\label{j_j-1Z2}
b_{j-1}^Z = b_j^Z - 1 \ge (r+1 - |I_{k+1}|) -1= \beta_{j}^Z - \epsilon_j \ge r+1 - |I_k| = \beta_{j-1}^Z,
\end{equation}
$$\beta_j^Y + \beta_j^Z - \epsilon_j = |I_k| + (r+1 - |I_{k+1}|) -1 \ge r+1.$$
Hence, (C2)--(C4) are satisfied.

Therefore, we conclude that $(\beta^Y, \beta^Z, \epsilon)$ is admissible. We claim this is what we are looking for, i.e., $\D(\beta^Y, \beta^Z, \epsilon) = \sum_{i=0}^r (a_i^Y + a_{r-i}^Z - d).$
By Proposition \ref{mainstr}, it is enough to show the following two conditions hold:
\begin{ilist}
\itm $J'_\beta$ is the empty set, or equivalently, $\beta_j^Y + \beta_{j-1}^Z = r+1,$ for any $j = 1, \dots, d.$
\itm For any $j \in J_\beta$ and any $i \in \psi_\beta^{-1}(j,j),$ we have that 
$(a^Y, a^Z)$ is connected at $i$ via $j,$
and $\epsilon_j =1.$
\end{ilist}
We have already seen that (i), i.e., \eqref{r+1cond} holds. Furthermore, by Corollary \ref{preimg}, \eqref{r+1cond}, \eqref{j_j-1Z1} and \eqref{j_j-1Z2},  we have that for any $j \in \{1, \dots, d-1\}$, $j \in J_\beta$ if and only if $\beta_j^Y + \beta_j^Z > r+1$ if and only if $\beta_j^Z > \beta_{j-1}^Z$ if and only if $j \in J.$ Thus,
$$J_\beta = J.$$

Suppose $j \in J_\beta = J$ and $i \in \psi_\beta^{-1}(j,j).$ Then $j = j_k,$  for some $k \in \{1, \dots, s\}.$ Therefore, by Lemma \ref{psibeta}, we have that
$$r+1-|I_k| = r+1 - \beta_j^Y  \le i < \beta_j^Z = r+1 - |I_{k+1}|.$$
However, by \eqref{j(i)Ik}, $r+1-|I_k|  \le i$ if and only if $\j(i) \ge j_k$, and $i <  r+1 - |I_{k+1}|$ if and only if $\j(i) < j_{k+1}.$ Hence, $j_k \le \j(i) < j_{k+1}.$ Since $\j(i) \in \nHat{J} \subseteq \{j_0 < j_1 < \dots < j_s < j_{s+1}\},$ $\j(i) = j_{k'}$ for a unique $k' \in \{0, \dots, s+1\}.$ Therefore, $\j(i) = j_k.$ Therefore, $(a^Y, a^Z)$ is connected at $i$ via $\j(i) = j_k = j.$ 
Moreover, by the definition of $\epsilon,$ it is clear that $\epsilon_j=1.$ Therefore, (ii) holds as well.
\end{proof}

Theorem \ref{main2} follows from Corollary \ref{nece} and Proposition \ref{suff}. Therefore, as we discussed in Section 2, we can conclude Theorem \ref{main}. We finish our paper with an example of a typical situation in which the conditions in Theorem \ref{main} are achieved.

\begin{ex}
Let $r =1$ and $d = 2.$ Set
$$a_0^Y = a_0^Z = 1, a_1^Y = a_1^Z = 2.$$
One checks that $(a^Y, a^Z)$ satisfies (A1)--(A3). Furthermore, $(a^Y, a^Z)$ is connected at both $i=0$ and $i=1$ via $j = 1.$ Thus, $(a^Y, a^Z)$ is connected. 

In the case when $g=0,$ there exist linear series with arbitrary vanishing sequences at a single fixed point, so in the reducible case, there exist EH limit linear series with arbitrary $(a^Y, a^Z)$ satisfying (A1)--(A3). In particular, we see that $G_d^r(X; a^Y, a^Z)$ is nonempty in the above example. Note that equality in (A3) does not hold for $i=0,1.$ Thus, $G_d^r(X; a^Y, a^Z)$ consists entirely of crude EH limit linear series. Hence, by Theorem \ref{main}, this gives examples of crude limit linear series, $\pi^{-1}(G_d^r(X; a^Y, a^Z))$,  which contain an open subset of $G_d^r(X)$.

\end{ex}

\subsection*{Acknowledgements}
I would like to thank Brian Osserman for suggesting this problem, and for explaining the geometric background.

\bibliographystyle{amsplain}
\bibliography{gen}

\end{document}